\documentclass{article}
\usepackage {amssymb}
\usepackage {amsfonts}
\usepackage {amsmath}
\usepackage {amsthm}
\usepackage{graphicx}
\usepackage{graphics}
\usepackage {framed}
\usepackage {amsxtra}
\usepackage {enumerate}
\usepackage[margin=1in]{geometry}

\newcommand{\restrict}{\mbox{$\mid\hspace{-1.1mm}\grave{}$}}

\theoremstyle{plain}
\newtheorem{thm}{Theorem}
\newtheorem{cor}[thm]{Corollary}

\theoremstyle{definition}

\theoremstyle{remark}

\begin{document}

\title{On sums of Egyptian fractions}

\author {Donald Silberger\footnote{SilbergD@newpaltz,edu\ \ State University of New York\ \ New Paltz, NY 12561}}

\maketitle

\[\mbox{\sf For our friend since 1948, Frank Oliver Wyse 1930-2024} \]

\begin{abstract} 

Let $n,d$, and $k$ be positive integers where $n$ and $d$ are coprime. Our two main results are

Theorem \ref{finites}. 
{\sl There is a partition  of the infinite interval $[kd,\infty)$ of positive integers into a family of finite sets $X$ for which the sum of the reciprocals of the elements in $X$ is $n/d$. }

Theorem \ref{infinites}. {\sl There is a partition of $[2kd,\infty)$ into an infinite family of infinite sets $Y$ for which the sum of the reciprocals of the elements in $Y$ is $n/d$.}

Our method is grounded in the {\sc vital identity}, $1/z = 1/(z+1) + 1/z(z+1)$, which holds for every complex number $z \notin \{-1,0\}$, and which gives rise to an eponymous algorithm that serves as our tool. 

At the core of our Theorems \ref{finites} and \ref{infinites}  is the number theoretic function $\star: x\mapsto \star x := x(x+1)$ into whose properties this paper continues an investigation initiated in \cite{Silbergers}.

\end{abstract} 

\section{Introduction and preliminaries} ``Egyptian fraction" is an ancient expression meaning  ``reciprocal $1/v$ of a positive integer $v$."

At least a century of interest in segments of the harmonic series precedes this paper. Our bibliography offers a brief history. The present work is a sequal to \cite{Silbergers}.  Our encompassing interests are the functions  

\[\sigma:{\cal P}({\mathbf N})\rightarrow{\mathbf Q}^+\cup\{\infty\}\quad\mbox{defined by}\quad \sigma: X\mapsto \sigma X := \sum_{x\in X} \frac{1}{x}\quad\mbox{and}\]
\[\Sigma:{\cal P}({\mathbf Q}^+)\rightarrow{\mathbf Q}^+\cup\{\infty\}\quad\mbox{defined by}\quad 
\Sigma: A\mapsto \Sigma A := \sum_{q\in A} q  \] where ${\mathbf N} := \{1,2,3,\ldots\}$, where ${\mathbf Q}^+$ is the set of positive rational numbers, and where  ${\tt P}(S)$ denotes the family of subsets of a given set $S$.  When $
{\cal F}$ is a family of sets, $\sigma{\cal F}$ denotes $\{\sigma X: X\in{\cal F}\}$.

Both functions were introduced in this paper's ``mother''\cite{Silbergers}, but our focus here is on $\sigma$. \vspace{.4em}

The methodology we employ deals with multisets $M$; some terminology thereof will be helpful.  We write $y^{(x)}\underline{\in} M$ to state that $y$ occurs in $M$ with {\it multiplicity} $x$; that is to say, $y$ occurs $x$ times as an element in $M$.  We will refer to $y^{(x)}$ as a ``$y$-tower that is $x$ stories, $y$, high.'' The expression $y^{(0)} \underline{\in} M$ asserts that $y\notin M$. The meaning of $\sigma$ expands in two natural ways: \  $\sigma y := \sigma \{y\} = 1/y$ and $\sigma y^{(x)}  := \sigma \{y^{(x)}\} := x\cdot\sigma \{y\} = x/y.$ Our multiset analog of set-algebra's union, $\cup$, is ``stack,'' $\sqcup$, whose meaning two examples should suffice to explain:  $\{1,2,3\}\sqcup\{2,3,4\} = \{1,2^{(2)},3^{(2)},4\}$ and $\{7^{(3)},8,9^{(5)},10\}\sqcup \{7,8,9^{(2)}\} = \{7^{(4)},8^{(2)},9^{(7)},10\}$.\vspace{.5em}

$[x,\infty)$ denotes $\{i : i\ge x\}\cap {\mathbf N}$ whereas ${\cal I}$ denotes the family of finite intervals $[i,i+j] := \{i,i+1,\ldots,i+j\}$ of consecutive positive integers. In 1915 
Theisinger\cite{Theisinger} proved that $\sigma[1,j]\in{\mathbf N}$ if and only if $j=1$. In 1918  K\"{u}rsch\'{a}k\cite{Kurschak} generalized this by proving that  ${\mathbf N}\cap\sigma{\cal I} = \{1\}$. Thus we learn that the function $\sigma\restrict{\cal I}$ is not surjective onto ${\mathbf Q}^+$. K\"{u}rsch\'{a}k's theorem is recalled as Exercise 3 on Page 7  in \cite{Baker}.  Other observations in this genre have been made by P. Erd\"{o}s\cite{Erdos2}, by H. Belbachir and A.  Khelladi\cite{Belbachir}, and by R. Obl\'{a}th\cite{Erdos2}. Each such result specifies an infinite family of finite sets $S$ of positive rationals for which $\Sigma S \notin{\mathbf N}$.

In 1946 Erd\"{o}s and Niven\cite{Erdos3}, veered toward another direction by proving that the function $\sigma\restrict{\cal I}$ is injective; i.e., that $\sigma[m,n] = \sigma[m',n']$ only if $m=m'$ and $n=n'$.

Two germane areas of study that suggest themselves are injectivity and surjectivity. We emphasize surgectivity here, and  improve on Theorem 1.1 in \cite{Silbergers}, which states that {\sl for every positive rational number $r$ there exists an infinite pairwise disjoint family of finite sets $S\subset{\mathbf N}$ for which $\sigma S = r$. } \vspace{.5em}

We use ideosyncratic terminology. The expression $[x,y]$ usually denotes the {\sf set} $\{z: x\le z \le y\}\cap{\mathbf N}$ but  sometimes, depending on context, it denotes instead the {\sf ascending sequence} $\langle x, x+1,\ldots,y\rangle$. Naturally, $[x,\infty)$ denotes the set (or the ascending sequence) of all integers $z\ge x$. 

Our  results depend upon the {\sc vital identity} $1/z = 1/(z+1)+1/z(z+1)$, which holds for every complex number $z \notin\{-1,0\}$. Introduced  in \cite{Silbergers} it was powered by two elementary number-theoretic functions ${\mathbf N}\rightarrow{\mathbf N}$; to wit:  the successor function $\diamond:i\mapsto i+1$ and a ``star" function $\star:i\mapsto i(i+1)$. So, the Vital Identity can be written as  $1/z = 1/\diamond z + 1/\star z$ for $z\in{\mathbf N}$. Despite its simplicity and elementary nature, the function $\star$ possesses  surprising  properties. We revisit $\star$ for its own sake later.

$\{\diamond,\star\}^j$ denotes the set of all length-$j$ words ${\sf w }$ in the alphabet $\{\diamond,\star\}$. Such a word ${\sf w}$ is treated also as the function created by the composition of the function-letters spelling the word ${\sf w}$. For a set $A\subset{\mathbf N}$ the expression ${\sf w}A$ denotes $\{{\sf w} x: x\in A\}$. We write $\{\diamond,\star\}^j A$ to mean $\{{\sf w}A: {\sf w}\in\{\diamond,\star\}^j\}$. Finally $\sigma \emptyset := 0$, and ${\sf w}\emptyset := \emptyset$ when ${\sf w}\in\{\diamond,\star\}^j$.

Both of our main theorems employ the functions $\diamond$ and $\star$ to produce a sequence of multiset modifications $M_i\rightarrow M_{i+1}$ via the replacement of a strategically selected integer $x\in M_i$ by the two integers $\diamond x := x+1$ and $\star x := x(x+1)$ since providentially $\sigma x = \sigma\diamond x+\sigma\star x$. Notice incidentally that $\diamond 1 = 1+1 = 1(1+1) = \star 1$, but fortunately $\diamond x = \star x\Leftrightarrow x=1$. \vspace{.3em}

Our proof of Theorem \ref{finites} uses the {\it Vital Algorithm} $\heartsuit$, defined below, to construct an infinite sequence of {\it transitional multisets}, the terms of a subsequence thereof being the member sets of the partition whose existence Theorem 1 asserts. The gears of $\heartsuit$ are the functions $\diamond$ and $\star$. Our proof of Theorem \ref{infinites} relies on insights gleaned from the proof of Theorem \ref{finites}.

\section{Main results}

\begin{thm}\label{finites} Let $\{d,n,k\}\subset{\mathbf N}$ with $d$ and $n$ coprime. Then there is a partition of the infinite set $[kd,\infty)$ of integers into finite sets $X$ for which $\sigma X = n/d$.
\end{thm}

\begin{proof} 

The Vital Algorithm, $\heartsuit := \{ \heartsuit_i^t: i,t$  nonnegative integers$\}$,  is a collection of multiset-valued functions $\heartsuit_i^t: M\mapsto M\heartsuit_i^t$  of a multiset variable $M$ whose operative rules are prescribed below. We will compose (only) these functions $\heartsuit_i^j$ from left to right; e.g., $M \heartsuit_i^u\circ \heartsuit_i^v = (M\heartsuit_i^u)\heartsuit_i^v$. Also, $S_i\heartsuit_i^j\bigcirc_{j=0}^t := S_i\heartsuit_i^0\circ \heartsuit_i^1\circ\cdots\circ\heartsuit_i^t$.

Now let $S_0 := \{kd\}\heartsuit_0^0 := \{kd\}$. Plainly, $\sigma S_0 = 1/kd$. For an integer $i\ge 0$ suppose that finite subsets $S_0,S_1,\ldots,S_i$ of $[kd,\infty)$ have been selected possessing the three following properties:\vspace{.3em} 

$\sigma S_j = 1/kd$ for every $j\in \{0,1,\ldots, i\}$.

The family $\{S_j: 0\le j\le i\}$ is pairwise disjoint.

$\bigcup_{j=0}^i S_j$ contains an initial segment at least $i$ integers long of the infinite sequence $[kd,\infty)$.  
\vspace{.4em}

Define $S_i\heartsuit_i^0 := S_i$ 

Define $S_{i+1} := S_i\heartsuit_i^u\bigcirc_{u=0}^{t_i}$ if  $t_i$ is the smallest positive integer $t$ for which $S_i\heartsuit_i^u\bigcirc_{u=0}^t$ is a simple set that satisfies  $S_i\heartsuit_i^u\bigcirc_{u=0}^t\cap \bigcup_{j=0}^i S_j = \emptyset$, and in this event define $S_{i+1}\heartsuit_{i+1}^0 := S_{i+1}$.

If either the multiset $S_i\heartsuit_i^u\bigcirc_{u=0}^t$ is not simple or $S_i\heartsuit_i^u\bigcirc_{u=0}^t\cap \bigcup_{j=0}^i S_j \not= \emptyset$, then two cases arise. Both cases pertain to the smallest integer $y$ for which $y^{(v)}\underline{\in}S_i\heartsuit_i^u\bigcirc_{u=0}^t$ for some $v\ge1$. Here we take it also that  $\{(y+1)^{(p)},(\star y)^{(q)}\}\subseteq S_i \heartsuit_i^u\bigcirc_{u=0}^t$ for some nonnegative integers $p$ and $q$.

In both cases, define $S_i\heartsuit_i^u\bigcirc_{u=0}^{t+1} := (S_i\heartsuit_i^u\bigcirc_{u=0}^t\setminus \{y^{(v)}\})\cup \{y^{(v-1)}, (y+1)^{(p+1)},(\star y)^{(q+1)}\}$. 

\underline{Case One}: $y\in\bigcup_{j=0}^i S_j$. Here $v=1$ is allowed and $y^{(0)}\underline{\in}S_i\heartsuit_i^u\bigcirc_{u=0}^{t+1}$; which says  $y\notin S_i\heartsuit_i^u\bigcirc_{u=0}^{t+1}$. So, ultimately $y\notin S_{i+1}$.  Case-One  maneuvers ensure that the family $\{S_j: 0\le j \le i+1\}$ is pairwise disjoint. 

\underline{Case Two}: $y\notin\bigcup_{j=0}^i S_j$. Here we require that $v>1$.  Case-Two maneuvers ensure that $y$ will be retained to become eventually an  element in $S_{i+1}$.
\vspace{.4em}

By the inductive hypothesis it is clear that $x\le y$ for every tower $y^{(x)} \underline{\in} S_i\heartsuit_i^u\bigcirc_{u=0}^t$. Moreover  $\star y - y = y^2$ for each positive integer $y$.  So, $x \le y = \star y -y^2 < \star y - y \le \star y - x$ for every $y>1$, and therefore a path of length $x$ that begins at $y$ must fail to reach $\star y$. Also, it is easily seen that $\star z < \star y-y^2$ whenever $z<y$.

Now take it  that $y\in \bigcup_{j=0}^i S_j$ and recall that $y$ is the smallest integer for which $y^{(v)} \underline{\in}  S_i\heartsuit_i^u\bigcirc_{u=0}^t$ with $v\ge1$. Case One pertains; $v=1$ may happen, whence $v(v-1)/2$ applications of the Vital Algorithm obliterates $y^{(v)}$  but the  $v(v-1)/2$ Algorithmic $\diamond$-actions result in $v$ towers built on the sites $y+1,y+2,\ldots, y+v$ along with a dispersed sequence $(\star y)^{(v)}, (\star(y+1))^{(v-1)}, \ldots, (\star(y+v-2))^{(2)},\star(y+v-1)$ of  new towers where no towers had existed before..  

If, instead, Case Two applies - i.e., if $y \notin\bigcup_{j=0}^i S_i$ - then  $v\ge 2$ is required. This time the Vital Algorithm reduces $y^{(v)}$ to $y^{(1)} = y$. After  $(v-1)(v-2)/2$ $\diamond$-actions by the Vital Algorithm starting with its assault on $y^{(v)}$,  the consequences are the single-story $v$ and $v-1$ towers on the sites $y,y+1,y+2,\ldots,y+v-1$. In addition, $v-1$ new towers  $(\star y)^{(v-1)}, (\star(y+1))^{(v-2)},\ldots,\star(y+v-2)$ will have arisen. 

In neither of these two cases is a newly created $\star$-built tower augmented by a tower built by or augmented by $\diamond$ actions.
 Thus we see that the maximum heights of towers occurring in the multisets $S_i\heartsuit_i^u\bigcirc_{u=0}^t$ is a nonincreasing function of $t$ and that it decreases  to $1$ for a sufficiently large  value $t_i$ of $t$. The multiset reached when $t=t_i$ is a simple finite set and is disjoint from $\bigcup_{j=1}^i S_j$. That is, $S_i\heartsuit_i^u\bigcirc_{u=0}^{t_i}= S_{i+1} =: S_{i+1}\heartsuit_{i+1}^0$.

The Vital Identity guarantees that $1/kd = \sigma S_i \heartsuit_i^u\bigcirc_{u=0}^t$ for all $t\ge 0$ by the inductive hypothesis.  Therefore 
$\sigma S_{i+1} =1/kd$. Furthermore, the family $\{S_j: 0\le j \le i+1\}$ is pairwise disjoint  The actions of the function $\diamond$ ensure that the maximum initial segment in $S_{i+1}$ of the infinite sequence $[kd,\infty)$  is a proper extention of the corresponding maximal initial segment thereof in $S_i$. It follows by induction that the Vital Algorithm produces an infinite partition $\{S_i: i= 0,1,2,\ldots\}$ of the set $[kd,\infty)$ into finite sets $S_i$ for each of which we have that $\sigma S_i = 1/kd$.

Finally, dice the infinite partition $\{S_0,S_1,S_2,\ldots\}$ into an infinite tribe $\{{\bf T}_u: u \in {\mathbf N}\}$ of $kn$-membered families ${\bf T}_u := \{S_{u_1},S_{u_2},\ldots, S_{u_{kn}}\}$. Let $X_u := \bigcup{\bf T}_u$ for each $u\in{\mathbf N}$. Then $\sigma X_u = kn/kd = n/d$ for every $u\in{\mathbf N}$.  The family $\{X_u : u \in {\mathbf N}\}$ is the partition of $[kd,\infty)$ whose existence the theorem asserts. \end{proof}

The following example illuminates our proof, exhibiting an ostensible snag and illustrating its resolution. \vspace{.3em}

\noindent{\bf Example One.} Let $k=d=n=1$. Since the Algorithm begins with $S_0 := \{kd\}$, in this example $S_0 = \{1\}$ and $\sigma S_0 = 1/kd =1$.  Replace the element $kd$ with the elements $\diamond (kd)$ and $\star (kd)$ thus transforming the set $S_0$ into the multiset $\{\diamond (kd),\star (kd)\}$, which happens not to be a simple set since $kd=1$. Instead,  $\{kd\}\heartsuit_0^1 = \{2^{(2)}\}$ because $\diamond (kd) =\diamond 1 = 2 = \star 1 =\star (kd)$.  In this example we can use more graceful terminology than the proof needed.  E.g., we write the foregoing multiset transformation in the form $S_0 := \{1\} \rightarrow_{0_1} \{2^{(2)}\}$.

The nonsimple multisets encountered in our construction process are transitional entities, stepping stones from an $S_i$ to $S_{i+1}$. In this example, we replace one of the two occurrences of $2$ in $\{2^{(2)}\}$ with the two elements $\diamond 2 = 2+1 =3$ and $\star 2 = 2(2+1) = 6$ thus creating the  simple set $S_1  := \{2,3,6\}$ having these three desired properties: $\sigma S_1  = 1$ and $\{1,2,3\}\subseteq S_0\cup S_1$ and $S_1\cap S_0 = \emptyset$. 

We now know  two, $S_0$ and $S_1$, of the infinitude of partition members $S_i$ we seek with $\sigma S_i = 1/kd = 1$ and such that  $S_0\cup S_1$ contains a length-2 initial segment of the sequence $[kd,\infty) = [1,\infty)$.

We modify $S_1$ until we obtain a simple set $S_2$ with the properties that $\sigma S_2 = 1$, that $S_2 \cap(S_0\cup S_1) = \emptyset$, and that $S_0\cup S_1\cup S_2$ contains a length-3 initial segment of $[kd,\infty)$. 

{\bf N.B.}  \  Always we choose to replace the \underline{smallest} integer from among those whose towers may come under assault in a multiset en route from an $S_i$ toward the construction of $S_{i+1}$. Thus in the present example, where $S_1 = \{2,3,6\}$, we replace $2$ with $\diamond 2 = 3$ and $\star 2 = 6$. This creates the transitional multiset $\{3^{(2)},6^{(2)}\}$. \vspace{.3em}

Summarizing and continuing: \ $S_0 = \{1\} \rightarrow_{0_1} \{2^{(2)}\} \rightarrow_{0_2} \{2,3,6\} =: S_1 \rightarrow_{1_1} \{3^{(2)}, 6^{(2)}\} \rightarrow_{1_2}$    \[   \{3,4,6^{(2)},12\} \rightarrow_{1_3}  \{4^{(2)},6^{(2)},12^{(2)}\} \rightarrow_{1_4} \{4,5,6^{(2)},12^{(2)},20\} \rightarrow_{1_5} \{4,5,6,7,12^{(2)},20,42\} \rightarrow_{1_6}     \]  \[ \{4,5,7^{(2)},12^{(2)},20,42^{(2)}\} \rightarrow_{1_7} \{4,5,7,8,12^{(2)},20,42^{(2)},56\} \rightarrow_{1_8} \{4,5,7,8,12,13,20, 42^{(2)},56.156\} \rightarrow_{1_9} \]  \[ \{4,5,7,8,12,13,20,42,43,56, 156,1806 =: S_2 \rightarrow_{2_1} \{5^{(2)},7,8,12,13,20^{(2)},  \,42,43,56,156,1806\} \rightarrow_{2_2}\ldots \]

In order to get from $S_1$ to $S_2$ we needed to ensure that $(S_0\cup S_1)\cap S_2 = \emptyset$ because we intend the infinite family  ${\cal F} := \{S_i: i=0,1,2,\ldots\}$ that is our penultimate target to be a partition of $[kd,\infty) = [1,\infty)$. So, ${\cal F}$ needs to be pairwise disjoint and therefore Case One  of the Vital Algorithm is activated when the Algorithmic chopper assaults a  $y^{(x)}$ for which $y$ is an element in an already-completed $S_i$. Case One guards against a previously incorporated  $y$ getting included also among the elements of the set $S_{i+1}$ that is  under construction. In the present example's construction of $S_0, S_1$, and $S_2$, the effect of Case One is manifested in the steps $\rightarrow_{0_1}, \rightarrow_{1_1}, \rightarrow_{1_3}$, and $\rightarrow_{1_6}$ and when headed from $S_2$ toward $S_3$ it occurs in step $\rightarrow_{2_1}$.

However, in order for ${\cal F}$ to be a partition of $[kd,\infty)$, it must also exhaust $[kd,\infty)$; i.e., it must satisfy $\bigcup{\cal F} = [kd,\infty)$. Every integer $y\ge kd$ must occur in some member of ${\cal F}$. Case Two guarantees this inclusivity. For, if $y$ does not occur as an element in any $S_j$ with $j\le i$ and the Algorithm's chopper attacks the tower $y^{(x)}$ with $x>0$ headed from $S_i$ toward  $S_{i+1}$, then the chopper does not destroy that tower but removes from it  exactly $x-1$ stories, thus preserving $y$ to be an element in the $S_{i+1}$ under construction from $S_i$. Eventually, for every integer $y\ge kd$ a tower $y^{(x)}$ with $x>0$ comes under assault by the Algorithm's chopper, which snuffs that tower's $y$ where appropriate but retains the $y$ when retention is appropriarte instead. 

The Algorithm's steps leading from $S_2$ to $S_3$ number in the hundreds. We have computed through to step $\rightarrow_{2_{109}}$. At that point our calculator was past its capacity to provide an exact integer product for $3274290\cdot3274291$ and our quest for $S_3$ ended. If someone with better electronics and more patience than ours wishes to complete the trip to $S_3$, then upon request we will mail or pdf those half dozen pages that we felt it prudent to omit from this paper. 

We are best reached at DonaldSilberger@gmail.com or at Sylvia.Silberger@gmail.com in the event that the author will have died. The author  is nearly ninety-five years old. It may be foolish to seek any $S_i$ with $i>5$, because as functions of $i$ we estimate that $\max S_i$ increases faster than $2^{2^i}$ and $\max |S_i|$ increases eventually faster than $2^i$. It has taken us until step $\rightarrow_{2_{75}}$ to demolish every tower $y^{(x)}$ for which  $y\in S_0\cup S_1\cup S_2$.
\vspace{.3em}

Now for a pedestrian, typical, less immediately troubled application to which Theorem \ref{finites} lends itself:

\vspace{.5em}

\noindent{\bf Example Two.} Let $k=1, n=2,$ and $d=3$. Then $S_0 = \{3\} \rightarrow_{0_1} \{4,12\} =: S_1 \rightarrow_{1_1} \{5,12,20\} \rightarrow_{1_2} \{5, 13, 20, 156\} =: S_2 \rightarrow_{2_1} \{6,13,20,30,156\} \rightarrow_{2_2} \{6,14,20,30,156,182\} \rightarrow_{2_3} \{6,14,21,30,156,182,420\}  \rightarrow_{2_4} \{6,14,21,30,157,182,420,24492\} =: S_3$. It will require $8$ steps to get from $S_3$ to $S_4$ and, after towers begin to form as they will when $i$ is large enough, then more than $2^i$ steps will be required from $S_i$ to $S_{i+1}$. 

One can verify that $\sigma S_i = 1/3$ for each of the four sets $S_i$ we have obtained in this example. Applying the Vital Agorithm repeatedly to the generation of these sets $S_i$, with each step requiring only half the time that its predecessor required, one acquires  after a finite period an infinite partition $\{S_0, S_1, S_2, \ldots \}$ of $[3,\infty)$ with $\sigma S_i = 1/3$ for each $i$. Then ${\cal G} :=\{S_0\cup S_1, S_2\cup S_3, S_4\cup S_5, \ldots \}$ is the ultimate partition of $[3,\infty)$ we seek here since obviously $\sigma(S_i\cup S_{i+1}) = 2/3$ for every $S_i\cup S_{i+1} \in {\cal G}$\vspace{.5em}

Thomas C. Brown asks whether, for each rational $r>0$, there is a partition of some subset of ${\mathbf N}$ into infinitely many infinite sets $Y$ for which $\sigma Y = r$? 

\begin{thm}\label{infinites} Let $k,n$, and $d$ be positive integers with $n$ and $d$ coprime. Then there is an infinite partition of $[2kd,\infty)$ into infinite sets $Y$ for which $\sigma Y = n/d$.
\end{thm}

One may wish to let $S_0 := \{2^jkd: j\in {\mathbf N}\}$ whence $\sigma S_0 = \sum_{j=1}^\infty 1/(2^jkd)= (1/kd)\sum_{j=1}^\infty 1/2^j = 1/kd$, leading to $S_1 = \bigcup_{j=1}^\infty\{2^jkd+1, 2^jkd(2^jkd+1)\}, \, S_2 = \bigcup_{j=1}^\infty\{\diamond^2 2^jkd, \diamond\star 2^jkd, \star\diamond 2^jkd, \star^2 2^jkd\}$, and so forth.  But troubles loom if we define every $S_i$ in this naive way: For large  $i$ we would wind up with multisets $S_i$ that are not simple. Moreover, the infinite family of all  $S_i$ thus obtained would fail to be pairwise disjoint. Our proof, below,  implements the idea behind this ingenuous approach while evading its pitfalls.

\begin{proof} Define $S_{0_1} := \{2kd\}$ and $S_{1_1} := \{ \diamond 2kd,\star 2kd\}$. Notice that $\sigma S_{0_1} = \sigma S_{1_1} = 1/2kd$. Clearly  $\{S_{0_1},S_{1_1}\}$ is a pairwise disjoint family of simple sets and $[2kd,2kd+1] := \{2kd,2kd+1\} \subseteq S_{0_1}\cup S_{1_1}$.

Define $S_{0_2} := \{2kd,4kd\}$ and $S^+_{1_2} := \{\diamond,\star\}^1S_{0_2}$ and $S^+_{2_2} :=\{\diamond,\star\}^2S_{0_2}$. If $\{S_{0_2},S^+_{1_2}\}$ is a pairwise disjoint family of simple sets then define $S_{1_2} := S^+_{1_2}$. Otherwise, apply the Vital Agorithm to $S^+_{1_2}$ a minimum number of times to achieve a simple set, $S$, such that $\{S_{0_2},S\}$ is pairwise disjoint, and then  rename $S_{1_2} := S$.  The minimality of Vital-Algorithm steps used here  ensures that $[2kd,2kd+1]\subseteq S_{0_2}\cup S_{1_2}$. Now ask whether $\{S_{0_2}, S_{1_2}, S^+_{2_2}\}$ is a pairwise disjoint family of simple sets; if ``yes" then rename $S_{2_2} :=S^+_{2_2}$.  If ``no,'' then Vital-Agorithmatize $S^+_{2_2}$ until achieving an $S$ that fits the bill, and then rename $S_{2_2} :=S$. Note that $[2kd,2kd+2]\subseteq S_{0_2}\cup S_{1_2}\cup S_{2_2}$ and also that $\sigma S_{u_2} = (1/2+1/4)/kd$ for all $u\in[0,2]$.

The preceeding two paragraphs provide the template for the  general step in our recursive construction.

Suppose for an  integer $i\ge2$ that an $i$-membered pairwise-disjoint family ${\cal S}_i := \{S_{0_i}, S_{1_i},\ldots, S_{i_i}\}$ of simple finite subsets $S_{v_i}$ has been constructed with the properties that: \[S_{0_i} = \{2^vkd: v\in[1,i]\}\,\,\mbox{and, for all}\,\, v\in[1,i],\,\,\,  \sigma S_{v_i} = \frac{1}{kd}\sum_{v=1}^i2^{-v}\, \mbox{and}\,\, [2kd,2kd+i] \subseteq \bigcup_{v=1}^iS_{v_i}. \]

Now let $S_{0_{i+1}} := S_{0_i}\cup\{2^{i+1}kd\}$ and for each $v\in[1,i+1]$ let $S^+_{v_{i+1}} := S_{v_i}\sqcup\{\diamond,\star\}^v\{2^{i+1}kd\}$. Finally, let $S^+_{{(i+1)}_{i+1}} := \{\diamond,\star\}^{i+1} S_{0_{i+1}}$. Suppose that the Vital Algorthm has been applied  - beginning with $u=1$, and proceeding from $u$ to $i+1$ for every $u\le i+1$ - in order to produce from $S^+_{u_{i+1}}$ a finite simple set $S_{u_{i+1}}$ with $2k+u\in S_{u_{i+1}}$, with $\sigma S_{u_{i+1}}= (kd)^{-1}\sum_{j=1}^{i+1} 2^{-j}$ and such that the family ${\cal S}_{i+1} :=\{S_{u_{i+1}}:u\in[0,i+1]\}$ is pairwise disjoint. Of course then $[2kd,2kd+i+1]\subseteq\bigcup_{v=1}^{i+1} S_{v_{i+1}}$.

The recursive construction in the two foregoing paragraphs provides for the creation of a pairwise disjoint family ${\cal S}_i = \{S_{0_i},S_{1_i},\ldots,S_{i_i}\}$ for every $i\in{\mathbf N}$ such that $[kd,kd+i]\subseteq\bigcup{\cal S}_i$ and such that $\sigma S_{0_i}=\sigma S_{u_i} = (1/kd)\sum_{j=1}^i (1/2^j)$ for all $u\in [0,i]$. As $i$ increases without bound the sizes of the families ${\cal S}_i$ likewise increase without bound. Moreover, when treating these sets $S_{u_i}$ as strictly increasing finite sequences we observe that the maximal eventually-constant prefixes of these $S_{u_i}$ increase  monotonically in length as $i$ increases so that ultimately the sets $S_{u_i}$ transmogrify as $i\rightarrow\infty$ into an infinite set $S_u$ for each $u$. The resulting infinite family ${\cal S} = \{S_0,S_1,S_2,\cdots\}$ is a partition of $[2kd,\infty)$ and for each $u$ we see that $\sigma S_u = (1/kd)\sum_{i=1}^\infty 1/2^i = 1/kd$.

We conclude our argument by partioning the family ${\cal S}$ into an infinite tribe of $kn$-membered subfamilies ${\sf T}_j$ whereupon the family $\{Y_j: j\in {\mathbf N}\}$ is the partition of $[2kd,\infty)$ of the sort whose existence the theorem asserts, when $Y_j := \bigcup {\sf T}_j$ for each $j$ since $\sigma Y_j = kn/kd = n/d$. \end{proof}\

\section{Mopping up}

There are two issues put forth in \cite{Silbergers} with which it is relevant for the present paper to deal.

\subsection{The function $\star$ }

If $x\mapsto x(x+1)$ were viewed as a real-valued function of a real variable we would have a parabola bringing to mind a comet from elsewhere in the Milky Way that visits our solar system once and once only. That function's  restriction, $\star:{\mathbf N}\rightarrow{\mathbf N}$, offers the lattice points on the first-quadrant portion of the  parabola.  

We focus is on the set $\{\star^p: p\in {\mathbf N}\}$ of positive-integer compositional powers of $\star$. Noticing that $\star x$ is a diophantine polynomial of degree $2$ we observe that $x\mapsto\star^px$ is a diophantine polynomial function of degree $2^p$ and that $n\mapsto\star ^nx$ is an exponential diophantine function of $n$.

\begin{thm}\label{star} For $x\in{\mathbf N}$ the integer $\star^nx$ has at least $n$ distinct prime factors. Moreover, $\star^nx|\star^mx$ for every $m>n$. Finally, if $p^k\| \star^nx$ for $p$ a prime, then $p^k\|\star^mx$ for every $m>n$.
\end{thm}

\begin{proof} Since $\star x = x(x+1)$ while the integers $x$ and $x+1$ are coprime and $x+1>1$, we see that $x+1$ has a prime factor that $x$ lacks. Therefore $\star^{n+1}x = \star^nx(\star^n x+1)$ has a prime factor that $\star^nx$ lacks.  It follows that $\star^n x$ has at least $n$ prime factors for every $n\in{\mathbf N}$. 

Clearly $\star^nx|\star^{n+1}x,\, \star^{n+1}x|\star^{n+2}x,\ldots,\star^{m-1}|\star^mx$. Thus when $n<m$ we get that $\star^nx|\star^mx$.

Let $p^k\|\star^nx$ with $p$ a prime.  Then $p$ is not a factor of any of the integers $\star^mx/\star^nx$ for which $m>n$.
\end{proof}

\begin{cor}\label{primes} The set of prime integers is infinite. 
\end{cor}

Theorem \ref{star} suggests an object that seems worthy of scrutiny; namely, the infinite $\langle\star,x\rangle$-sequence of prime powers that we define as follows: 

First, express the prime-powers factorization of $x$ in the order of increasing prime sizes. Then express the prime-powers factorization of $x+1$ in the order of ascending prime sizes. Ditto the prime power factorization of $\star x+1$ in the order of ascending prime sizes. Do the same for $\star^2 x+1$, next for $\star^3x+1$, then for $\star^4x+4$ and so on forever and concatenate one string of prime-powers factorizations in the order that they are obtained as you obtain them.  The resulting infinite concatenation of prime-powers factorizations is $x^\star$.\vspace{.5em} 

\noindent{\bf Example Three.} Let $x :=2$. We compute an initial segment (aka ``preface'') of $2^\star$, demarcating the disinct prime-power factorizations with the punctuation  \  {\bf\large ;}  \  \ Thus we begin $2^\star = 2;3;7;43;13\cdot139;3263443;\ldots$

Is it coincidence that four out of the five ``principal,'' $\star^i2+1$, segments are primes while the fifth is the square-free product of two primes? That seems to portend a high concentration of principal segments that are either themselves primes or that have unusually few factors. Of course, for small $x$ it would have to happen that many of the small primes or small powers thereof would occur as factors of the principal segments, and Theorem \ref{star} ensures that those used-up primes will not recur as factors of any subsequent principal segments. But already the fifth principal factor of $2^\star$ is a stand-alone prime larger then three million; that seems striking. The sequences $x^\star$ are rich sources of open questions. For instance:\vspace{.5em}

\noindent{\bf Questions One.} With ${\mathbf P}$ denoting the set of all prime integers and ${\mathbf Pp}$ as the set of all prime powers and ${\mathbf P}[x]$ and ${\mathbf Pp}[x]$ denoting respectively the set of prime factors and the set of prime-power factors of $x$ or of some $\star^ix+1$, what can be said about the sundry values of $x$ with respect to the latter two sets? E.g.:

For which values if any of $x$ does ${\mathbf P}[x] = {\mathbf P}$ occur?

For which values if any of $x\not=y$ does  ${\mathbf P}[x] ={\mathbf P}[y]$? 

For which values if any of $x\not= y$ does ${\mathbf Pp}[x] = {\mathbf Pp}[y]$? 

\subsection{Harmonic series revisited}

It is easily proved but also corollary to Khinchin's Theorem\cite{Khinchin} that the set ${\mathbf H}$ of all finite segment sums $\sigma[n,n+k]$ of the harmonic series is dense in the positive real numbers.\vspace{.5em}

\noindent{\bf Question Two.}  Is the set ${\mathbf Q}^+\setminus{\mathbf H}$ dense in the positive real numbers?

\vspace{2em}

\noindent{\bf 2020 Mathematics Subject Classification:} \, 01A55, 01A60, 11A25, 11N05, 11N25


\begin{thebibliography}{99}



\bibitem{Baker} A. Baker, ``A Concise Introduction to the Theory of Numbers'', Cambridge University Press (1984).

\bibitem{Belbachir} H. Belbachir and A. Khelladi, {\sl On a sum involving powers of reciprocals of an arithmetic progression}, Ann. Mathematicae et Informaticae {\bf 34} (2007), 29-31. 

\bibitem{Erdos1} P. Erd\"os, {\sl Egy K\"ursch\'ak-F\'ele Elemi Sz\'amelm\'eleti T\'etel \'Altal\'anos\'it\'asa}, Matematikai \'es Fizikai Lapok  BD. XXXIX, Budapest (1932), 1-8.  

\bibitem{Erdos2} P. Erd\"os, {\sl A theorem of Sylvester and Schur}, J. London Math. Soc. {\bf 9} (1934), 191-258. 

\bibitem{Erdos3} P. Erd\"{o}s and I. Niven, {\sl Some properties of partial sums of the harmonic series}, Bull. Amer. Math. Soc. {\bf 52} (1946), 248-251.

\bibitem{Euler} L. Euler, {\sl Variae observationes circa series infinites}, (1744), Theorem 19, {\it Euler Archive - All Works}, 72. 
https://scholarlycommons.pacific.edu/euler-works/72.

\bibitem{Silbergers} David Hobby, Donald Silberger, Sylvia Silberger {\sl Sums of finitely distinct rationals}, arXiv:1702.01316v2[math.NT] 19 Feb 2019, 1-12.  

\bibitem{Hoffman} P. Hoffman, ``The Man Who Loved Only Numbers: The Story of Paul Erd\"os and the Search for Mathematical Truth'',  N. Y. Hyperion, 1998.

\bibitem{Khinchin} A. Y. Khinchin, {\sl Einige S\"{a}tze \"{u}ber die Kettenbr\"{u}che mit Anwendungen auf die Theorie der Dipphantischen Approximationen,} Mathematische Annalen {\bf 92}(1-2) (1924), 115-125.


\bibitem{Kurschak} J. K\"ursch\'ak, Matematikai \'es Fizikai Lapok, {\bf 27} (1918), 299. 

\bibitem{Shorey} T. N. Shorey, {\sl Theorems of Sylvester and Schur}, Math. Student (2007), Special Centenary Volume (2008), 135-145. Online article (http://www.math.tifr.res.in/ shorey/newton.pdf).

\bibitem{Schur} I. Schur, {\sl  Einige S\"atze \"uber Primzahlen mit Anwendungen auf Irreduzibilit\"atsfragen. II.} Sitzungsber. Preuss. Akad. Wiss. Berlin Phys. Math. K1. {\bf 14} (1929), 370-391.

\bibitem{Sylvester} J. J. Sylvester, {\sl On Arithmetical Series}, Messenger of Math. {\bf 21} (1892), 1-19, 87-120 (Collected Mathematical Papers, Bd. {\bf 14}, 687-731).  

\bibitem{Theisinger} L. Theisinger, {\sl Bemerkung \"uber die harmonische Reihe}, Monatshefte f\"ur Mathematik und Physik {\bf 26} (1915), 132-134. 

\end{thebibliography}
\end{document}